\documentclass[a4paper,11pt]{article}

\usepackage{amsmath,amssymb}

\usepackage{graphicx}
\usepackage{epsfig}
\usepackage{color}

\newenvironment{proof}{\begin{trivlist}
	\item[\noindent]{\it Proof\:}}{\quad $\square$\end{trivlist}}

\newenvironment{exl}{\begin{trivlist}
	\item[\noindent]{\bf Example\:}}{\end{trivlist}}

\newtheorem{dfn}{Definition}[section]
\newtheorem{thm}{Theorem}
\newtheorem{prp}{Proposition}
\newtheorem{lem}{Lemma}

\def\R{{\mathbb R}}

\def\phi{\varphi}
\def\epsilon{\varepsilon}

\def\E{{\mathcal E}}
\def\F{{\mathcal F}}
\def\C{{\mathcal C}}

\def\grad{\mathrm{grad}\,}
\def\ext{\mathrm{ext}}

\title{A variational proof of Alexandrov's\\ convex cap theorem}
\author{Ivan Izmestiev
\thanks{Research for this article was supported by the DFG Research Unit 565 ``Polyhedral Surfaces''.}\\
Institut f\"ur Mathematik\\
Technische Universit\"at Berlin\\
Str. des 17. Juni 136\\
10623 Berlin, Germany\\
{\tt izmestiev@math.tu-berlin.de}}

\begin{document}

\maketitle

\begin{abstract}
We give a variational proof of the existence and uniqueness of a convex cap with the given upper boundary. The proof uses the concavity of the total scalar curvature functional on the space of generalized convex caps. As a byproduct, we prove that generalized convex caps with the fixed boundary are globally rigid, that is uniquely determined by their curvatures.
\end{abstract}

\section{Introduction}
In this paper we present a new proof of the following theorem by A.~D. Alexandrov \cite[Section 5.1, Theorem 2]{Al05}.
\begin{thm} \label{thm:Alex}
Let $D$ be a disk with a convex Euclidean polyhedral metric. Then there exists a convex cap $C \subset \R^3$ with the upper boundary isometric to $D$. Besides, $C$ is unique up to a rigid motion.
\end{thm}

Informally speaking, a {\it convex cap} is a convex polytope that lies on a plane and does not throw a shadow when the sun is in the zenith. The upper boundary $D$ of a convex cap is topologically a disk. The metric on $D$ is locally Euclidean except at the vertices. Indeed, a point on an edge has a small neighborhood isometric to a Euclidean disk, the two halves of the disk lying on the two faces adjacent to the edge. Similarly, neigborhoods of the interior vertices are isometric to Euclidean cones with angles less than $2\pi$. The boundary of $D$ is locally modelled on the half-plane and angular regions with angles less than $\pi$. Any metric that has these properties is called a {\it convex Euclidean polyhedral metric}. Thus the upper boundary of any convex cap is a disk with a convex Euclidean polyhedral metric. Theorem \ref{thm:Alex} is the converse: any metric of this kind can be realized as the upper boundary of a convex cap, and this is in an essentially unique way.

\subsection{Idea of the proof and organization of the paper}
In Section \ref{sec:Prelim}, we give formal definitions and make some preliminary remarks.

Our proof uses the \textit{variational method}. In this method the object whose existence is to be proved is identified with a critical point of a functional. If the functional achieves its maximum in the interior of the domain, then the maximum point is a critical point and yields a desired object. If, moreover, the functional is concave and the domain convex, then there is only one critical point. This implies the uniqueness of the object.

The domain of our functional is the space $\C(D)$ of \textit{generalized convex caps} with the upper boundary $D$; it is studied in Section \ref{sec:GCC}. The functional on $\C(D)$ is the \textit{total scalar curvature} $S$ which is studied in Section \ref{sec:S}. Section \ref{sec:Proof} contains the details of the proof.

Subsection \ref{subsec:GCCDef} introduces generalized convex caps. Roughly speaking, a generalized convex cap is a convex cap that has cone singularities along the heights. Heights are the segments that join the vertices of the upper boundary with the base and are orthogonal to the base. If the curvatures around heights vanish, then the generalized convex cap becomes a usual convex cap. To define a generalized cap, one needs to specify a triangulation $T$ of $D$ and a collection $h$ of heights. However, in Subsection \ref{subsec:HeightsToCaps} we show that the heights already suffice. This makes the space $\C(D)$ naturally a subset of $\R^\Sigma$, where $\Sigma$ is the singular set of $D$. In Subsection \ref{subsec:SpaceOfCaps} we show that $\C(D)$ is in fact a convex polytope in $\R^\Sigma$.

In Subsection \ref{subsec:TSC} we define the total scalar curvature of a generalized convex cap. It is similar to the definition
$$
\sum_e \ell_e(\pi - \theta_e)
$$
of the total mean curvature of a convex polytope, where $\ell_e$ is the length of the edge $e$ and $\theta_e$ is the angle at $e$. For generalized caps one takes additionally the sum of the heights multiplied with their curvatures. In Subsection \ref{subsec:SConc} we prove that the function $S$ is concave and in important cases strictly concave. This allows us to prove two rigidity theorems \ref{thm:GenRig} and \ref{thm:InfRigCap} for generalized and classical convex caps in the next Subsection. Also, due to the Darboux-Sauer theorem this yields a new proof of the infinitesimal rigidity of convex polytopes.

The most part of Subsection \ref{subsec:Proof} deals with the case when $S$ achieves its maximum on the boundary of $\C(D)$. This happens when some faces of the resulting convex cap are vertical. In Subsection \ref{subsec:Volkov} we sketch the proof of Theorem \ref{thm:Alex} by Volkov from \cite{Vol60}.

\subsection{Related work} \label{subsec:RelWork}
The variational method was extensively used in the study of circle packings and patterns, see \cite{BS2} and references therein. The functionals that appear there can often be interpreted as volumes of certain hyperbolic polyhedra.

Our primary goal was to find a variational proof of Alexandrov's theorem on the existence and uniqueness of a polytope with the given metric on the boundary, see Subsection \ref{subsec:AlVol}. A similar approach, via generalized convex polytopes and the total scalar curvature meets with an obstacle there: the functional is neither concave nor convex. In return, a rather unexpected relation arises between the total scalar curvature and the volume of the dual. This allows one to prove the non-degeneracy of the Hessian of the total scalar curvature by extending the Alexandrov-Fenchel inequalities for mixed volumes. This variation on the variational method leads to a proof presented in \cite{BI}.

Alexandrov's theorems are, in fact, special cases of a very general statement about realization of a Euclidean, spherical, or hyperbolic metric with cone singularities of positive or negative curvature on a surface of an arbitrary genus. Each time, such a metric can be realized in a unique way as the boundary of a certain polyhedral object. The case of a sphere with the singularities of positive curvature is Alexandrov's theorem on Euclidean, hyperbolic and spherical polytopes. Spherical metrics on the sphere with singularities of negative curvature were considered in \cite{HR}. For the formulation of the general statement and the proof of the case of genus greater than one, see \cite{Fil}. The case of metrics on the torus is the subject of the paper \cite{FI}.

\subsection{Proofs by Alexandrov and Volkov} \label{subsec:AlVol}
Theorem \ref{thm:Alex} looks similar to another, more renowned theorem of Alexandrov: \textit{let $S$ be a sphere with a convex Euclidean polyhedral metric; then there exists a unique convex polytope $P \subset \R^3$ with the boundary isometric to $S$.} In his book, Alexandrov derives the cap theorem from the polytope theorem. Indeed, by identifying two copies of $D$ along the boundary one gets a convex polyhedral metric on the sphere. Because of the uniqueness, the polytope that realizes this metric must have a plane of symmetry. Thus it can be cut in two convex caps. This proves the existence of a cap with the upper boundary $D$. The uniqueness follows from the uniqueness of the polytope with the boundary the doubled $D$.

To prove his polytope theorem, Alexandrov used a method called sometimes the \textit{deformation method} or the invariance of domain method. Consider the space ${\mathcal P}$ of all convex polytopes and the space ${\mathcal M}$ of all convex polyhedral metrics on the sphere. There is an obvious map ${\mathcal P} \to {\mathcal M}$ that associates to a polytope the metric of its boundary. The existence and uniqueness of a polytope with a given metric is translated as the bijectivity of this map. Alexandrov shows that the map is a local homeomorphism, which is related to the infinitesimal rigidity. Then he proves certain topological properties of the spaces and of the map which imply that the map is in fact a bijection.

Alexandrov applied the deformation method to several existence and uniqueness statements for polyhedral objects, see \cite{Al05}. Recently it was used in \cite{HR}, \cite{Fil} to prove parts of the general statement mentioned at the end of Subsection \ref{subsec:RelWork}. A drawback of this method is that it is not constructive: it provides no algorithm how to construct the polytope with a given development. Alexandrov points this out in \cite[pp. 320--321]{Al05} after he proves by deformation Minkowski's theorem on the existence and uniqueness of a convex polytope with given face normals and face areas. Alexandrov notes that Minkowski's own proof is by the variational method, and poses the problem to find a similar proof for his polytope theorem.

Yu.~A. Volkov, a student of Alexandrov, found constructive proofs both for the polytope and for the cap theorem. He considered generalized convex caps with positive curvatures around the heights and showed that the curvatures vanish when the sum of heights is maximized. This proof is published in \cite{Vol60}, for English translation see \cite[Section 12.2]{Al05}. Volkov's proof of the polytope theorem is similar and is the subject of Volkov's PhD of 1955, see also \cite{VP}. Volkov's approach is not variational in the above sense. The functional is linear, so the maximum point always lies on the boundary of the domain and is not a critical point of the functional. But his proof is elementary and, in the case of the cap, quite short.

Our proof is to a large extent inspired by the works of Volkov. In \cite{Vol68}, reprinted as \cite[Section 12.1]{Al05} Volkov studies the local behavior of the function $S$ and proves claims that are very close to the concavity of $S$. Completely new in our paper are the explicit description of the space $\C(D)$ as a compact convex subset of a Euclidean space and the global rigidity statement for generalized convex caps.

\subsection{Acknowledgements}
I am indebted to Alexander I. Bobenko for bringing the problem to my attention and for helpful discussions. I thank Fran\c{c}ois Fillastre and Jean-Marc Schlenker for useful remarks.

\section{Preliminaries} \label{sec:Prelim}
Here we give the formal definitions of what was sketched in the Introduction and reduce Theorem \ref{thm:Alex} to the case when the polyhedral disk $D$ has at least one interior and at least one boundary singularity.

\begin{dfn} \label{dfn:CC}
A {\rm convex cap} is a convex polytope $C$ in $\R^3$ with the following properties:
\begin{enumerate}
\item $C$ is contained in the upper half-space $\R^3_+$, and $C \cap \partial \R_+^3 \ne \emptyset$. The face $C \cap \partial \R^3$ of $C$ is called the base of the cap $C$;
\item the orthogonal projection $\R^3 \to \R^2 = \partial \R_+^3$ maps $C$ to its base;
\item the dimension of the polytope $C$ is 2 or 3.
\end{enumerate}
\end{dfn}

\begin{dfn}
Let $C$ be a convex cap. If $\dim C = 3$, then the {\rm upper boundary} of $C$ is $\partial C$ minus the interior of the base. If $\dim C = 2$ and $C \subset \partial \R_+^3$, then the upper boundary of $C$ is $C$. Finally, if $\dim C = 2$ and the base of $C$ is an edge of the polygon $C$, then the upper boundary consists of two copies of $C$ identified along the edges except the base edge.
\end{dfn}

\begin{figure}[ht]
\begin{center}
\begin{picture}(0,0)%
\includegraphics{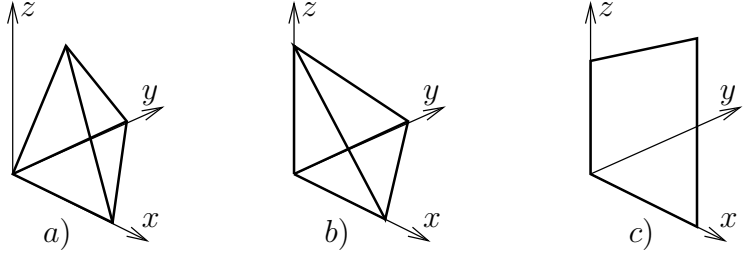}%
\end{picture}%
\setlength{\unitlength}{4144sp}%
\begingroup\makeatletter\ifx\SetFigFont\undefined%
\gdef\SetFigFont#1#2#3#4#5{%
  \reset@font\fontsize{#1}{#2pt}%
  \fontfamily{#3}\fontseries{#4}\fontshape{#5}%
  \selectfont}%
\fi\endgroup%
\begin{picture}(4503,1503)(589,-1186)
\put(1396,-1051){\makebox(0,0)[lb]{\smash{{\SetFigFont{12}{14.4}{\rmdefault}{\mddefault}{\updefault}{\color[rgb]{0,0,0}$x$}%
}}}}
\put(3061,-1051){\makebox(0,0)[lb]{\smash{{\SetFigFont{12}{14.4}{\rmdefault}{\mddefault}{\updefault}{\color[rgb]{0,0,0}$x$}%
}}}}
\put(4816,-1051){\makebox(0,0)[lb]{\smash{{\SetFigFont{12}{14.4}{\rmdefault}{\mddefault}{\updefault}{\color[rgb]{0,0,0}$x$}%
}}}}
\put(4276,-1141){\makebox(0,0)[lb]{\smash{{\SetFigFont{12}{14.4}{\rmdefault}{\mddefault}{\updefault}{\color[rgb]{0,0,0}$c)$}%
}}}}
\put(676,209){\makebox(0,0)[lb]{\smash{{\SetFigFont{12}{14.4}{\rmdefault}{\mddefault}{\updefault}{\color[rgb]{0,0,0}$z$}%
}}}}
\put(2341,209){\makebox(0,0)[lb]{\smash{{\SetFigFont{12}{14.4}{\rmdefault}{\mddefault}{\updefault}{\color[rgb]{0,0,0}$z$}%
}}}}
\put(4096,209){\makebox(0,0)[lb]{\smash{{\SetFigFont{12}{14.4}{\rmdefault}{\mddefault}{\updefault}{\color[rgb]{0,0,0}$z$}%
}}}}
\put(1396,-286){\makebox(0,0)[lb]{\smash{{\SetFigFont{12}{14.4}{\rmdefault}{\mddefault}{\updefault}{\color[rgb]{0,0,0}$y$}%
}}}}
\put(3061,-286){\makebox(0,0)[lb]{\smash{{\SetFigFont{12}{14.4}{\rmdefault}{\mddefault}{\updefault}{\color[rgb]{0,0,0}$y$}%
}}}}
\put(4816,-286){\makebox(0,0)[lb]{\smash{{\SetFigFont{12}{14.4}{\rmdefault}{\mddefault}{\updefault}{\color[rgb]{0,0,0}$y$}%
}}}}
\put(811,-1141){\makebox(0,0)[lb]{\smash{{\SetFigFont{12}{14.4}{\rmdefault}{\mddefault}{\updefault}{\color[rgb]{0,0,0}$a)$}%
}}}}
\put(2476,-1141){\makebox(0,0)[lb]{\smash{{\SetFigFont{12}{14.4}{\rmdefault}{\mddefault}{\updefault}{\color[rgb]{0,0,0}$b)$}%
}}}}
\end{picture}%
\end{center}
\caption{The upper boundary of the cap b) has one interior and two boundary singularities; the cap c) is a polygon, its upper boundary has two interior and no boundary singularities.}
\label{fig:ExCaps}
\end{figure}

Figure \ref{fig:ExCaps} shows several examples of convex caps.

\begin{dfn} \label{dfn:EPM}
Let $M$ be a surface, possibly with boundary. A metric structure on $M$ is called a {\rm Euclidean polyhedral metric} if there is a finite set $\Sigma \subset M$ of points called singular points, such that
\begin{enumerate}
\item
any regular interior point $x$ has a neighborhood isometric to an open subset of the Euclidean plane; any regular boundary point $x$ has a neigborhood isometric to an open subset of the half-plane;
\item
any singular interior point $x$ has a neighborhood isometric to an open subset of a cone, with $x$ at the apex of the cone; any singular boundary point $x$ has a neighborhood isometric to an open subset of an angular region, with $x$ at the angle's vertex.
\end{enumerate}
A Euclidean polyhedral metric is called {\rm convex}, if all of the angles at the interior singularities are less that $2\pi$, and all of the angles at the boundary singularities are less that $\pi$.
\end{dfn}
For brevity, we often omit the word ``Euclidean'' when speaking of polyhedral metrics.

Let $D$ be the upper boundary of a convex cap $C$. It is not hard to show that $D$ is a disk with a convex polyhedral metric. Any vertex of $C$ that does not lie in the base is an interior singularity of $D$. For a vertex $v$ of $C$ that lies in the base there are two possibilities: if $v$ belongs to a vertical edge of $C$, then $v$ is regular, otherwise it is singular.

\begin{lem}
Theorem \ref{thm:Alex} is true when $D$ has no interior singularities. The case where $D$ has no boundary singularities follows from the case where there is at least one boundary singularity.
\end{lem}
\begin{proof}. If $D$ has no interior singularities, then $D$ is isometric to a convex polygon. Clearly, any convex cap with the upper boundary $D$ is obtained by isometrically embedding $D$ in $\partial \R_+^3$.

Let $D$ have no boundary singularities. It follows that the $\epsilon$-neighborhood of $\partial D$ is isometric to $\partial D \times [0,\epsilon)$ for some $\epsilon > 0$. Take the maximum $\epsilon$ with this property. The locally isometric embedding $\partial D \times [0,\epsilon) \to D$ extends to a map
$$
\iota: \partial D \times [0,\epsilon] \to D.
$$
Assume that the map $\iota$ is injective. Then the complement
$$
D' = D \setminus \iota(\partial D \times [0,\epsilon))
$$
is a disk with a convex polyhedral metric and the boundary $\iota(\partial D \times \{\epsilon\})$. Due to the choice of $\epsilon$, $\partial D'$ contains a singularity. By Theorem \ref{thm:Alex} applied to the disk $D'$, there is a convex cap $C'$ with the upper boundary $D'$. Lifting the cap $C'$ to the height $\epsilon$, we obtain a convex cap $C$ with the upper boundary $D$. Conversely, if $C$ is any convex cap with the upper boundary $D$, then $C \cap (\partial \R_+^3 \times [0,\epsilon])$ is a right prism. Therefore the uniqueness of $C$ follows from the uniqueness of $C'$.

Assume that the map $\iota$ is not injective. The convexity of $D$ implies that the map $\iota$ folds the circle $\partial D \times \{\epsilon\}$, so that $D$ is isometric to two copies of a rectangle identified along three sides. The corresponding convex cap is this rectangle standing on its fourth side. The uniqueness is not hard to prove. \end{proof}

\begin{dfn}
For a set $A \subset \R^3_+$, its \emph{lower hull} is the convex hull of the union $A \cup \mathrm{pr}(A)$, where $\mathrm{pr}: \R^3_+ \to \R^2 = \partial \R^3_+$ is the orthogonal projection.
\end{dfn}
From the discussion after Definition \ref{dfn:EPM} it follows that a convex cap is the lower hull of the singularities of its upper boundary.

\section{Generalized convex caps} \label{sec:GCC}

\subsection{Definitions} \label{subsec:GCCDef}
Let $D$ be a convex polyhedral disk such that $\Sigma \cap \partial D \ne \emptyset$, where $\Sigma$ is the singular set of $D$. We denote the elements of $\Sigma$ by the letters $i,j,\ldots$. A \emph{geodesic triangulation} $T$ of $D$ is a decomposition of $D$ into triangles by geodesics with the endpoints in $\Sigma$. By $\E(T)$ and $\F(T)$ we denote the sets of edges and triangles of $T$, respectively.

Note that we allow multiple edges as well as loops. In particular, triangles of $T$ may have identifications on the boundary. An edge with endpoints $i$ and $j$ is denoted by $ij$, and a triangle with vertices $i, j, k$ by $ijk$. This notation may be ambigous, but it should not lead to confusion.

An edge $ij$ of $T$ is called a \emph{boundary edge}, if it is contained in the boundary of $D$; otherwise it is called an \emph{interior edge}.

\begin{exl} Consider the convex cap which is the lower hull of the points $(1,0,0), (0,1,0)$, and $(0,0,1)$, see Figure \ref{fig:ExCaps} b). Its upper boundary $D$ has one interior singularity $(0,0,1)$ and two boundary singularities $(1,0,0)$ and $(0,1,0)$. Disk $D$ has a unique geodesic triangulation. It consists of two triangles: an equilateral one with side length $\sqrt{2}$ and an isosceles right one with leg length $\sqrt{2}$. The triangles have two edges in common.
\end{exl}

Any convex polyhedral disk can be geodesically triangulated. This can be proved in several different ways. For example, draw the shortest geodesics from a boundary singularity to all interior singularities. After cutting along them we have a (possibly non-convex, possibly self-overlapping) polygon. It is a classical result that every polygon can be triangulated by diagonals.

\begin{prp} \label{prp:FinTriang}
The number of geodesic triangulations of $D$ is finite.
\end{prp}
\begin{proof}. By \cite[Section 4.5]{BZ}, the lengths of geodesics in $D$ are bounded from above. On the other hand, the lengths of geodesics joining two singular points of $D$ form a discrete set. It follows that the set of geodesics between singular points is finite. Therefore, the number of geodesic triangulations is also finite. \end{proof}

\begin{dfn} \label{dfn:Prism}
A {\rm prism} is a convex polytope isometric to the lower hull of three non-collinear points in $\R^3_+$.
\end{dfn}
We use the term ``prism'' not in its usual meaning: the lateral edges of our prism are necessarily orthogonal to the lower base, but the upper base need not be parallel to it. Up to isometry, a prism is uniquely determined by a Euclidean triangle $ijk$ isometric to its upper base, and by three non-negative \emph{heights} $h_i, h_j, h_k$, which are the lengths of the lateral edges. Note that a height may equal 0.

\begin{figure}[ht]
\begin{center}
\begin{picture}(0,0)%
\includegraphics{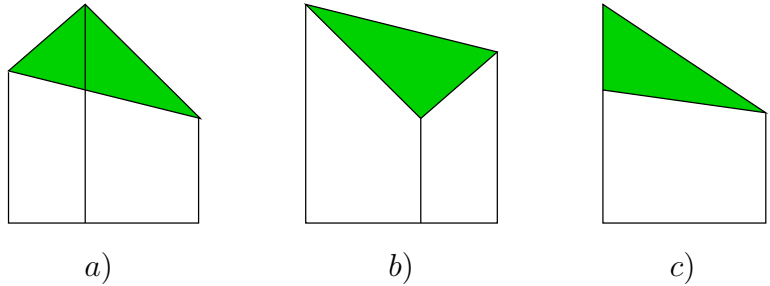}%
\end{picture}%
\setlength{\unitlength}{4144sp}%
\begingroup\makeatletter\ifx\SetFigFont\undefined%
\gdef\SetFigFont#1#2#3#4#5{%
  \reset@font\fontsize{#1}{#2pt}%
  \fontfamily{#3}\fontseries{#4}\fontshape{#5}%
  \selectfont}%
\fi\endgroup%
\begin{picture}(4506,1677)(439,-1411)
\put(901,-1366){\makebox(0,0)[lb]{\smash{{\SetFigFont{12}{14.4}{\rmdefault}{\mddefault}{\updefault}{\color[rgb]{0,0,0}$a)$}%
}}}}
\put(2701,-1366){\makebox(0,0)[lb]{\smash{{\SetFigFont{12}{14.4}{\rmdefault}{\mddefault}{\updefault}{\color[rgb]{0,0,0}$b)$}%
}}}}
\put(4366,-1366){\makebox(0,0)[lb]{\smash{{\SetFigFont{12}{14.4}{\rmdefault}{\mddefault}{\updefault}{\color[rgb]{0,0,0}$c)$}%
}}}}
\end{picture}%
\end{center}
\caption{Three types of degenerate prisms.}
\label{fig:DegPrisms}
\end{figure}

If the upper base of a prism is orthogonal to the lower base, the prism degenerates to a polygon. In this polygon we still distinguish the triangle that is the upper base of the prism. A degeneration can occur in three different ways shown on Figure \ref{fig:DegPrisms}. In the cases a) and b) we can assign values $0$ or $\pi$ to the dihedral angles at the lateral edges and at the edges of the upper base. In the case c) some of these angles are not defined.

\begin{dfn} \label{dfn:GCC}
A {\rm generalized convex cap} $C$ with the upper boundary $D$ is a polyhedron glued from prisms (in the sense of Definition \ref{dfn:Prism}) whose upper bases are the triangles of a geodesic triangulation $T$ of $D$. The identification pattern of the prisms corresponds to the combinatorics of the triangulation. Besides, the following properties should hold:
\begin{enumerate}
\item the heights of the boundary vertices are 0, i.~e. in a prism that contains a boundary edge $ij$ this edge is shared by the upper and the lower base;
\item for every interior edge $ij \in \E(T)$, the dihedral angle $\theta_{ij}$ is either not defined or does not exceed $\pi$. Here $\theta_{ij}$ is the sum of the two dihedral angles of the prisms at the edge $ij$.
\end{enumerate}
\end{dfn}
At the left of Figure \ref{fig:Prisms} two prisms are shown that share an interior edge~$ij$.

Lemma \ref{lem:NoVertEdg} in the next subsection ensures that in a generalized convex cap the angles $\theta_{ij}$ are always defined.

Clearly, a generalized convex cap can be given by a couple $(T,h)$, where $T$ is a geodesic triangulation of $D$, and $h:\Sigma \to \R$ is an assignment of heights to the vertices of $T$. A couple $(T,h)$ produces a convex cap if and only if the following conditions are satisfied. First, $h_i = 0$ if $i \in \partial D$. Second, for any triangle $ijk \in \F(T)$ there are restrictions on $h_i, h_j, h_k$ ensuring that a prism with the upper base $ijk$ and heights $h_i, h_j, h_k$ exist. And third, the angles $\theta_{ij}$ must be less than or equal to $\pi$, if defined.

Two different couples $(T,h)$ and $(T', h')$ can produce the same generalized convex cap. For example, if $\theta_{ij} = \pi$ for some edge $ij \in \E(T)$, and the edge $ij$ is shared by two triangles that form a convex quadrilateral, then the triangulation $T$ can be changed by replacing the edge $ij$ through another diagonal of the quadrilateral. On the other hand, $(T,h) \sim (T',h')$ implies $h = h'$ and $ij \in \E(T')$ for all $ij \in \E(T)$ such that $\theta_{ij} < \pi$. Therefore every generalized convex cap defines a tesselation of the disk $D$ obtained by erasing in any associated triangulation $T$ the edges $ij$ with $\theta_{ij} = \pi$. In Subsection \ref{subsec:HeightsToCaps} we show that the collection of heigths $h$ defines the cap uniquely.

\begin{exl} Take any triangulation $T$ and put $h_i = 0$ for all $i$. The result is a generalized convex cap with $\theta_{ij} = \pi$ for all $ij$. \end{exl}

It is easy to see that the following quantity is well-defined.

\begin{dfn} \label{dfn:kappa}
Let $(T,h)$ be a generalized convex cap. For any interior singularity $i \in \Sigma \setminus \partial D$ denote by $\omega_i$ the sum of the dihedral angles of the prisms at the edge under the vertex $i$. The angle defect
$$
\kappa_i = 2\pi - \omega_i
$$
is called the {\rm curvature} at the $i$-th height.
\end{dfn}
In the above example with all of $h_i = 0$ the curvature $\kappa_i$ equals the angle defect of the singularity $i$ in the metric of $D$.

We refer to a convex cap in the sense of Definition \ref{dfn:CC} as a classical convex cap. 

Let $C$ be a classical convex cap with at least one singularity on the boundary. If $C$ has no vertical faces, it is easy to represent it as a generalized convex cap: it suffices to triangulate non-triangular faces of $C$ and then cut $C$ into non-degenerate prisms. If $C$ has vertical faces, then we first remove all of the vertical edges, if there are any. The remaining edges subdivide $D$ into polygons with vertices in $\Sigma$. Since any polygon can be triangulated by diagonals, we can complete this subdivision to a geodesic triangulation $T$. Thus we obtain a set of prisms where some prisms are degenerate. If $\dim C = 3$, then the cap $C$ is isometric to the union of the non-degenerate prisms, and the degenerate prisms cover the vertical faces of $C$ ``from outside''. If $\dim C = 2$, one has to distinguish two cases: when $C$ has no vertical edges, and when $C$ has one. In the latter case the triangulation $T$ contains a triangle ``bent across'' the vertical edge of $C$.

As an example, consider Figure \ref{fig:ExCaps} b). Here we have a non-degenerate prism and a degenerate one. The degenerate prism is formed by the two vertical faces of the cap. The curvature $\kappa_i$ at the unique interior singularity equals $\frac{\pi}{2}$.

\subsection{Heights define a cap} \label{subsec:HeightsToCaps}
Here and later on we use the word ``cap'' to refer to a generalized convex cap.

Let us denote
$$
\C(D) = \{\mbox{generalized convex caps with the upper boundary }D\}.
$$
The main result of this Subsection is Proposition \ref{prp:HeightsToCaps} which is implied by Lemmas \ref{lem:CapsAreFunctions} and \ref{lem:FuncOnSing}.

\begin{dfn}
A function $f:D \to \R$ is called a {\rm PL function} if it is linear on every triangle of some geodesic triangulation of $D$.

A function $f:D \to \R$ is called {\rm concave} if its restriction to any geodesic segment in $D$ is concave as the function of the arc length.
\end{dfn}

\begin{lem} \label{lem:CapsAreFunctions}
The space $\C(D)$ is in a one-to-one correspondence with the space of concave PL functions $\widetilde{h}:D \to \R$ such that
\begin{equation}
\widetilde{h}|_{\partial D} = 0 \label{eqn:hTilde=0}
\end{equation}
and
\begin{equation}
|\grad\widetilde{h}(x)| \le 1 \label{eqn:hGrad}
\end{equation}
at every point $x \in D$ such that $\grad\widetilde{h}(x)$ is defined.

If $C \in \C(D)$ is a cap represented by $(T,h)$, then the corresponding function $\widetilde{h}$ assigns to the point $x \in D$ the distance from $x$ to the lower base of a prism that contains $x$. In particular,
$$
\widetilde{h}(i) = h_i
$$
for every $i \in \Sigma$.
\end{lem}
\begin{proof}. Take any $C \in \C(D)$. Let $\widetilde{h}:D \to \R$ be the distance function from the upper boundary to the base of the cap $C$, as described in the proposition. It is easy to see that $\widetilde{h}$ is a concave PL function that satisfies conditions (\ref{eqn:hTilde=0}) and (\ref{eqn:hGrad}).

Conversely, let $\widetilde{h}:D \to \R$ be any concave PL function that satisfies (\ref{eqn:hTilde=0}) and (\ref{eqn:hGrad}). Concavity and property (\ref{eqn:hTilde=0}) imply that $\widetilde{h}$ is non-negative inside $D$. Let $T$ be a geodesic triangulation of $D$ such that $\widetilde{h}$ is linear on the triangles of $T$. Take any triangle $ijk \in \F(T)$. Property (\ref{eqn:hGrad}) implies that there exists a prism with the upper base $ijk$ such that $\widetilde{h}$ is the distance function from the upper to the lower base. If $ijl$ is another triangle in $T$, then the sum of the dihedral angles of the prisms at the edge $ij$, if defined, is less than or equal to $\pi$ due to the concavity of $\widetilde{h}$.

It is not hard to see that two caps with equal distance functions $\widetilde{h}$ are equal. \end{proof}

\begin{lem} \label{lem:NoVertEdg}
A cap cannot contain degenerate prisms of types b) and c) from Figure \ref{fig:DegPrisms}. In particular, the angles $\theta_{ij}$ and $\kappa_i$ are defined for any cap.
\end{lem}
\begin{proof}. Assume that the cap contains a prism of type c). Let $ij$ be the vertical edge of the upper base so that $i$ lies above $j$. Then $|\grad \widetilde{h}| = 1$ inside both triangles adjacent to the edge $ij$, and the direction of the gradient is parallel to the edge $ij$. Note that the vertex $j$ cannot lie on the boundary of $D$, because then it would not be a singular point. Consider two geodesics on both sides of the edge $ij$ and parallel to $ij$. Follow them in the direction of decrease of the function $\widetilde{h}$. If the geodesics are sufficiently close to the edge $ij$, then they intersect at a point $x \in D$. Then the directional derivatives of $\widetilde{h}$ in two different directions at $x$ are equal to 1. Therefore $|\grad \widetilde{h}(x)| > 1$ that contradicts (\ref{eqn:hGrad}).

The proof that a degenerate prism of type b) cannot occur is similar. \end{proof}

\begin{lem} \label{lem:FuncOnSing}
The extension of a map $h: \Sigma \to \R$ to a concave PL function $\widetilde{h}: D \to \R$ is unique, if exists.
\end{lem}
\begin{proof}. Assume that there are two different extensions $\widetilde{h}$ and $\widetilde{h}'$. Without loss of generality, there is a point $x \in D$ such that
$$
\widetilde{h}(x) > \widetilde{h}'(x).
$$
Let $T$ be a geodesic triangulation such that $\widetilde{h}$ is linear on the triangles of $T$. Let $ijk \in \F(T)$ be a triangle that contains $x$. Then $\widetilde{h}(x)$ is the linear interpolation from the values $h_i, h_j, h_k$ at the vertices $i, j, k$. Consider the function $\widetilde{h}'$ on the triangle $ijk$. At the vertices $i, j, k$ it takes the same values as the function $\widetilde{h}$. Besides, $\widetilde{h}'$ is concave. This implies $\widetilde{h}'(x) \ge \widetilde{h}(x)$ that contradicts our assumption. \end{proof}

\begin{prp} \label{prp:HeightsToCaps}
The map $(T,h) \mapsto h$ defines a bijection between the space $\C(D)$ and the space of maps $h: \Sigma \to \R$ such that
\begin{enumerate}
\item $h$ admits an extension to a concave PL function $\widetilde{h}: D \to \R$;
\item the map $\widetilde{h}$ has the properties (\ref{eqn:hTilde=0}) and (\ref{eqn:hGrad}).
\end{enumerate}
\end{prp}

\subsection{The space of generalized convex caps} \label{subsec:SpaceOfCaps}
By Proposition \ref{prp:HeightsToCaps}, the space $\C(D)$ is identified with a subset of the Euclidean space $\R^{\Sigma}$. Here we show that $\C(D)$ is a bounded convex polyhedron in $\R^{\Sigma}$.

\begin{figure}[ht]
\begin{center}
\begin{picture}(0,0)%
\includegraphics{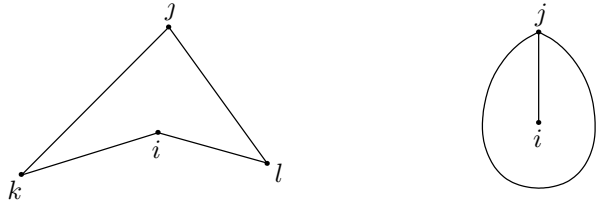}%
\end{picture}%
\setlength{\unitlength}{4144sp}%
\begingroup\makeatletter\ifx\SetFigFont\undefined%
\gdef\SetFigFont#1#2#3#4#5{%
  \reset@font\fontsize{#1}{#2pt}%
  \fontfamily{#3}\fontseries{#4}\fontshape{#5}%
  \selectfont}%
\fi\endgroup%
\begin{picture}(3506,1275)(5976,-1363)
\put(6924,-211){\makebox(0,0)[lb]{\smash{{\SetFigFont{10}{12.0}{\familydefault}{\mddefault}{\updefault}{\color[rgb]{0,0,0}$j$}%
}}}}
\put(7575,-1207){\makebox(0,0)[lb]{\smash{{\SetFigFont{10}{12.0}{\familydefault}{\mddefault}{\updefault}{\color[rgb]{0,0,0}$l$}%
}}}}
\put(6847,-1060){\makebox(0,0)[lb]{\smash{{\SetFigFont{10}{12.0}{\familydefault}{\mddefault}{\updefault}{\color[rgb]{0,0,0}$i$}%
}}}}
\put(5991,-1308){\makebox(0,0)[lb]{\smash{{\SetFigFont{10}{12.0}{\familydefault}{\mddefault}{\updefault}{\color[rgb]{0,0,0}$k$}%
}}}}
\put(9101,-994){\makebox(0,0)[lb]{\smash{{\SetFigFont{10}{12.0}{\familydefault}{\mddefault}{\updefault}{\color[rgb]{0,0,0}$i$}%
}}}}
\put(9115,-247){\makebox(0,0)[lb]{\smash{{\SetFigFont{10}{12.0}{\familydefault}{\mddefault}{\updefault}{\color[rgb]{0,0,0}$j$}%
}}}}
\end{picture}%
\end{center}
\caption{Every concave quadrilateral $ikjl$ gives rise to an inequality (\ref{eqn:ConcQuadr}), every triangle $iji$ gives rise to an inequality (\ref{eqn:CloGeod}).}
\label{fig:fig1}
\end{figure}

We start by studying a larger space $\C'(D)$ which is obtained from $\C(D)$ by ignoring conditions (\ref{eqn:hTilde=0}) and (\ref{eqn:hGrad}).
\begin{dfn}
The space $\C'(D) \subset \R^{\Sigma}$ consists of the maps $h: \Sigma \to \R$ that admit a concave PL extension to $D$.
\end{dfn}
To state the next Proposition, we need to introduce some notions. A Euclidean quadrilateral $ikjl$ in $D$ is a region that is bounded by simple geodesic segments $ik, kj, jl, li,$ and that contains no singularities in the interior. Vertices and segments are allowed to coincide. A Euclidean quadrilateral can be developed onto the plane, after resolving possible identifications on the boundary. A Euclidean triangle $jij$ arises when there is a geodesic loop based at $j$ that encloses a unique singularity $i$. See Figure \ref{fig:fig1}.

Define the function $\ext_{jkl}$ on the quadrilateral $ikjl$ as the linear function that takes values $h_j,\, h_k,\, h_l$ at the respective vertices.
\begin{lem} \label{lem:C'(D)}
The space $\C'(D)$ is a convex polyhedron in $\R^{\Sigma}$. Namely, it is the solution set of a system of linear inequalities of the form:
\begin{eqnarray}
h_i & \ge & \ext_{jkl}(i), \label{eqn:ConcQuadr}\\
h_i & \ge & h_j. \label{eqn:CloGeod}
\end{eqnarray}
There is one equation of the form (\ref{eqn:ConcQuadr}) for each Euclidean quadrilateral $ikjl$ with the angle at $i$ greater or equal $\pi$, and one equation of the form (\ref{eqn:CloGeod}) for each Euclidean triangle $jij$, see Figure \ref{fig:fig1}.
\end{lem}

\begin{proof}. Assume that $h$ admits a concave PL extension $\widetilde{h}$. Then the concavity of $\widetilde{h}$ easily implies inequalities (\ref{eqn:ConcQuadr}) and (\ref{eqn:CloGeod}).

Let us prove the sufficiency of conditions (\ref{eqn:ConcQuadr}) and (\ref{eqn:CloGeod}). Let $T$ be any geodesic triangulation. Denote by $\widetilde{h_T}$ the PL function that is linear on the triangles of $T$ and takes value $h_i$ at every $i \in \Sigma$. Call an edge of $T$ \emph{good}, if the function $\widetilde{h_T}$ is concave across this edge, otherwise call an edge \emph{bad}. Our aim is to find a triangulation with good edges only.

Apply the {\it flip algorithm} to the triangulation $T$. Let $ij$ be an interior edge of $T$. If it belongs to two different triangles $ijk$ and $ijl$, and if the quadrilateral $ikjl$ is strictly convex, then the edge $ij$ can be {\it flipped}. To flip means to replace the edge $ij$ through the diagonal $kl$ of the quadrilateral $ikjl$. In the flip algorithm we start with an arbitrary triangulation, pick up a bad edge, flip it, look for a bad edge in the new triangulation, flip it and so on. The flip algorithm terminates when there are no bad edges.

Conditions (\ref{eqn:ConcQuadr}) and (\ref{eqn:CloGeod}) ensure that if $ij$ is a bad edge, then it can be flipped. Let us show that the flip algorithm terminates. Note that when a bad edge is flipped, the function $\widetilde{h_T}$ increases pointwise. By Proposition \ref{prp:FinTriang}, the polyhedral disk $D$ has only finitely many geodesic triangulations. Therefore, the flip algorithm cannot run infinitely. \end{proof}

\begin{thm} \label{thm:C(D)}
The space $\C(D)$ is a non-empty bounded convex polyhedron in~$\R^{\Sigma}$. Namely, it is the set of points that satisfy conditions (\ref{eqn:ConcQuadr}), (\ref{eqn:CloGeod}), and
\begin{eqnarray}
h_i & = & 0 \ \quad \mbox{ for all }i \in \partial D, \label{eqn:hBound}\\
h_i & \le & d_i \quad \mbox{ for all }i, \label{eqn:hInt}
\end{eqnarray}
where $d_i$ is the distance in $D$ from $i$ to $\partial D$.
\end{thm}
\begin{proof}. The necessity of (\ref{eqn:hBound}) is obvious. To show the necessity of (\ref{eqn:hInt}), draw a shortest geodesic $ix$ joining $i$ to $\partial D$. The segment $ix$ has length $d_i$. The restriction of $\widetilde{h}$ to $ix$ is a PL function with the absolute value of derivative less than or equal to 1. Since $\widetilde{h}(x) = 0$, we have $h_i = \widetilde{h}(i) \le d_i$.

It remains to show that conditions (\ref{eqn:hTilde=0}) and (\ref{eqn:hGrad}) are fulfilled for any concave PL function $\widetilde{h}$ whose values at the points of $\Sigma$ satisfy (\ref{eqn:hBound}) and (\ref{eqn:hInt}). Property (\ref{eqn:hTilde=0}) follows from (\ref{eqn:hBound}) and the piecewise linearity of $\widetilde{h}$. To show (\ref{eqn:hGrad}), assume the converse: there is an $x \in D$ with $|\grad\widetilde{h}(x)| = c > 1$. Since the gradient is constant inside every triangle of $T$, we can choose the point $x$ arbitrarily close to some singular point $i$. Draw a geodesic that starts at $x$ and goes in the direction of the fastest decrease of $\widetilde{h}$ at $x$. This geodesic ends either at a boundary point $y$ of $D$ or at a singularity. By perturbing the point $x$ we can assume that the former is the case. Then the restriction of $\widetilde{h}$ to the segment $xy$ is a concave function that vanishes at $y$. Since its slope at $x$ equals $-c$, we have
$$
\widetilde{h}(x) \ge c \cdot |xy|.
$$
By taking the limit when $x$ tends to $i$, we obtain
$$
h_i \ge c \cdot d_i > d_i,
$$
which contradicts (\ref{eqn:hInt}).

The polyhedron $\C(D)$ is bounded because it is contained in the box $\prod_{i \in \Sigma} [0,d_i]$. The inequality $0 \le h_i$ for every $i$ follows from (\ref{eqn:hBound}) and from the concavity of the function $\widetilde{h}$.

Finally, $\C(D)$ is non-empty because it contains the point $h=0$. \end{proof}

There is a natural decomposition
\begin{equation} \label{eqn:Decomp}
\C(D) = \bigcup_T \C^T(D),
\end{equation}
where $\C^T(D)$ consists of the caps that have a representative of the form $(T,h)$.
\begin{prp}
For every geodesic triangulation $T$, the space $\C^T(D)$ is a bounded convex polyhedron in $\R^\Sigma$.
\end{prp}
\begin{proof}. Let $h: \Sigma \to \R$ be a point of $\C(D)$. It belongs to $\C^T(D)$ if and only if the PL extension $\widetilde{h_T}$ of $h$ with respect to $T$ is concave. Let $ij$ be an interior edge of $T$. If it belongs to two different triangles $ijk$ and $ijl$, then $\widetilde{h_T}$ is concave across the edge $ij$ if and only if
\begin{equation} \label{eqn:CT}
\ext_{ijk}(l) \ge h_l.
\end{equation}
If $ij$ is incident to one triangle only, then we have the situation on the right of Figure \ref{fig:fig1}, up to exchanging $i$ and $j$. In this case the concavity of $\widetilde{h_T}$ across $ij$ is equivalent to
\begin{equation} \label{eqn:CT1}
h_i \ge h_j.
\end{equation}

Therefore $\C^T(D)$ is the set of points $h \in \C(D)$ that satisfy the system of linear inequalities of the type (\ref{eqn:CT}) and (\ref{eqn:CT1}), one inequality for every edge of~$T$. \end{proof}

\section{Total scalar curvature and rigidity} \label{sec:S}

\subsection{The total scalar curvature of a generalized convex cap} \label{subsec:TSC}
\begin{figure}[ht]
\begin{center}
\begin{picture}(0,0)%
\includegraphics{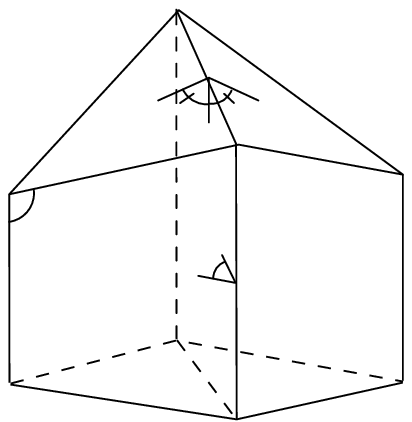}%
\end{picture}%
\setlength{\unitlength}{4144sp}%
\begingroup\makeatletter\ifx\SetFigFont\undefined%
\gdef\SetFigFont#1#2#3#4#5{%
  \reset@font\fontsize{#1}{#2pt}%
  \fontfamily{#3}\fontseries{#4}\fontshape{#5}%
  \selectfont}%
\fi\endgroup%
\begin{picture}(1992,2047)(1707,-1266)
\put(2912, 81){\makebox(0,0)[lb]{\smash{{\SetFigFont{10}{12.0}{\familydefault}{\mddefault}{\updefault}{\color[rgb]{0,0,0}$\alpha_{ji}$}%
}}}}
\put(2967,-675){\makebox(0,0)[lb]{\smash{{\SetFigFont{10}{12.0}{\familydefault}{\mddefault}{\updefault}{\color[rgb]{0,0,0}$\omega_{ijk}$}%
}}}}
\put(2427, 90){\makebox(0,0)[lb]{\smash{{\SetFigFont{10}{12.0}{\familydefault}{\mddefault}{\updefault}{\color[rgb]{0,0,0}$\alpha_{ij}$}%
}}}}
\put(1707,-315){\makebox(0,0)[lb]{\smash{{\SetFigFont{10}{12.0}{\familydefault}{\mddefault}{\updefault}{\color[rgb]{0,0,0}$k$}%
}}}}
\put(2022,-360){\makebox(0,0)[lb]{\smash{{\SetFigFont{10}{12.0}{\familydefault}{\mddefault}{\updefault}{\color[rgb]{0,0,0}$\rho_{ki}$}%
}}}}
\put(2535,649){\makebox(0,0)[lb]{\smash{{\SetFigFont{10}{12.0}{\familydefault}{\mddefault}{\updefault}{\color[rgb]{0,0,0}$j$}%
}}}}
\put(2962,-159){\makebox(0,0)[lb]{\smash{{\SetFigFont{10}{12.0}{\familydefault}{\mddefault}{\updefault}{\color[rgb]{0,0,0}$i$}%
}}}}
\end{picture}%
\end{center}
\caption{Angles and lengths in a generalized convex cap.}
\label{fig:Prisms}
\end{figure}

\begin{dfn}
Let $C$ be a generalized convex cap represented by $(T,h)$. The total scalar curvature of $C$ is defined as
$$
S(C) = \sum_{\Sigma \setminus \partial D} h_i \kappa_i + \sum_{\mathrm{int}D} \ell_{ij}(\pi - \theta_{ij}) + \sum_{\partial D} \ell_{ij}(\frac{\pi}{2} - \eta_{ij}).$$
Here $\kappa_i = 2\pi - \omega_i$ is the curvature at the $i$-th height of the cap $C$, $\ell_{ij}$ is the length of the edge $ij \in \E(T)$, $\theta_{ij}$ is the total dihedral angle at an interior edge $ij$, and $\eta_{ij}$ is the dihedral angle at a boundary edge $ij$. See Figure \ref{fig:Prisms}.
\end{dfn}
By Lemma \ref{lem:NoVertEdg}, the angles $\kappa_i$ and $\theta_{ij}$ are defined for any polyhedral complex $(T,h)$ that represents a generalized convex cap.

Let $(T',h)$ be another representative of the same cap $C$. If $ij$ is an edge of $T$ but not of $T'$, then we have $\theta_{ij} = \pi$ and thus $\ell_{ij}(\pi - \theta_{ij}) = 0$. The curvature $\kappa_i$ does not depend on the choice of a representative of $C$. Therefore the function $S$ is well-defined and continuous on $\C(D)$.

\begin{figure}[ht]
\begin{center}
\begin{picture}(0,0)%
\includegraphics{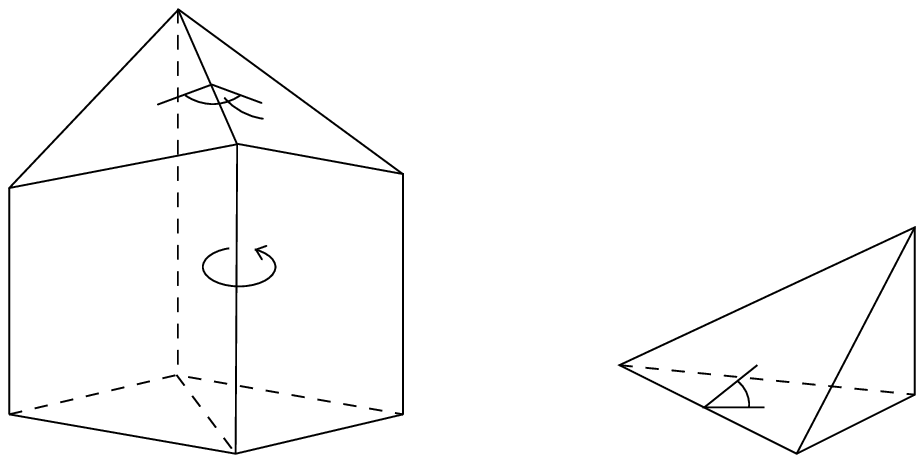}%
\end{picture}%
\setlength{\unitlength}{4144sp}%
\begingroup\makeatletter\ifx\SetFigFont\undefined%
\gdef\SetFigFont#1#2#3#4#5{%
  \reset@font\fontsize{#1}{#2pt}%
  \fontfamily{#3}\fontseries{#4}\fontshape{#5}%
  \selectfont}%
\fi\endgroup%
\begin{picture}(4391,2432)(482,-1604)
\put(1916, 55){\makebox(0,0)[lb]{\smash{{\SetFigFont{10}{12.0}{\familydefault}{\mddefault}{\updefault}{\color[rgb]{0,0,0}$\theta_{ij}$}%
}}}}
\put(1786,-184){\makebox(0,0)[lb]{\smash{{\SetFigFont{10}{12.0}{\familydefault}{\mddefault}{\updefault}{\color[rgb]{0,0,0}$i$}%
}}}}
\put(1414,696){\makebox(0,0)[lb]{\smash{{\SetFigFont{10}{12.0}{\familydefault}{\mddefault}{\updefault}{\color[rgb]{0,0,0}$j$}%
}}}}
\put(1981,-601){\makebox(0,0)[lb]{\smash{{\SetFigFont{10}{12.0}{\familydefault}{\mddefault}{\updefault}{\color[rgb]{0,0,0}$\omega_i$}%
}}}}
\put(904,270){\makebox(0,0)[lb]{\smash{{\SetFigFont{10}{12.0}{\familydefault}{\mddefault}{\updefault}{\color[rgb]{0,0,0}$\ell_{jk}$}%
}}}}
\put(535,-151){\makebox(0,0)[lb]{\smash{{\SetFigFont{10}{12.0}{\familydefault}{\mddefault}{\updefault}{\color[rgb]{0,0,0}$k$}%
}}}}
\put(482,-781){\makebox(0,0)[lb]{\smash{{\SetFigFont{10}{12.0}{\familydefault}{\mddefault}{\updefault}{\color[rgb]{0,0,0}$h_k$}%
}}}}
\put(4316,-1564){\makebox(0,0)[lb]{\smash{{\SetFigFont{12}{14.4}{\rmdefault}{\mddefault}{\updefault}{\color[rgb]{0,0,0}$i$}%
}}}}
\put(3358,-1021){\makebox(0,0)[lb]{\smash{{\SetFigFont{12}{14.4}{\rmdefault}{\mddefault}{\updefault}{\color[rgb]{0,0,0}$j$}%
}}}}
\put(4161,-993){\makebox(0,0)[lb]{\smash{{\SetFigFont{12}{14.4}{\rmdefault}{\mddefault}{\updefault}{\color[rgb]{0,0,0}$\eta_{ij}$}%
}}}}
\end{picture}%
\end{center}
\caption{More angles.}
\label{fig:Prisms1}
\end{figure}

\begin{dfn}
Let $C$ be a generalized convex cap represented by $(T,h)$. For any $i \ne j \in \Sigma$ put
$$
a_{ij} = \left\{ \begin{array}{ll}
		 \frac{\cot \alpha_{ij} + \cot \alpha_{ji}}{\ell_{ij} \sin^2 \rho_{ij}} & \mbox{if } ij \mbox{ is an interior edge of }T;\\
		 0 & \mbox{otherwise.}
		 \end{array}
\right.
$$
Here $\alpha_{ij}$ and $\alpha_{ji}$ are the dihedral angles of the prisms at the edge $ij$, thus $\alpha_{ij} + \alpha_{ji} = \theta_{ij}$; $\rho_{ij}$ is the angle between the edge $ij$ and the $i$-th height. See Figure \ref{fig:Prisms1}.

If $h_i = 0$, then the angle $\rho_{ij}$ is defined as the angle between the edge $ij$ and the vector $(0,0,-1)$ at the vertex $i$.

If there are several interior edges in $T$ that join $i$ and $j$, then $a_{ij}$ is the sum of the corresponding expressions over all such edges.
\end{dfn}

If $\theta_{ij} = \pi$, then $\cot \alpha_{ij} + \cot \alpha_{ji} = 0$. Therefore $a_{ij}$ does not depend on the choice of a triangulation $T$. Besides, due to $\rho_{ij} + \rho_{ji} = \pi$ we have $a_{ij} = a_{ji}$.

The quantities $a_{ij}$ will appear in the expressions for the second partial derivatives of the function $S$, see Proposition \ref{prp:DerS} below. Since the domain $\C(D)$ of $f$ is a polyhedron, the notion of a smooth function needs to be suitably generalized.

\begin{dfn} \label{dfn:DiffFunc}
Let $f: X \to \R$ be a continuous function on a polyhedron $X \subset \R^n$. We say that $f$ is of class $C^1$ on $X$ and write $f \in C^1(X)$, iff the directional derivatives $\frac{\partial f}{\partial \xi}$ exist for all admissible directions $\xi$, and depend on $\xi$ linearly.

More exactly, we require the existence of continuous functions $f_i: X \to \R$ for $i = 1,\ldots, n$ such that for any $x \in X$ and any $\xi \in \R^n$ such that $x + \epsilon \xi \in X$ for all sufficiently small positive $\epsilon$, we have
$$
\frac{\partial f}{\partial \xi} (x) = \sum_{i = 1}^n f_i \xi_i.
$$
We say that the function $f$ is of class $C^2$ iff $f \in C^1(X)$ and $f_i \in C^1(X)$ for all~$i$.
\end{dfn}
We call functions $f_i$ partial derivatives of $f$ even if they don't make sense as directional derivatives. This is the case exactly when $x + \epsilon e_i \notin X$ for all positive $\epsilon$, where $e_i$ denotes the $i$-th basis vector.

Definition \ref{dfn:DiffFunc} is equivalent to saying that $f \in C^1(X)$ iff
$$
f(x + dx) = f(x) + \sum_{i=1}^n f_i dx_i + o(\|dx\|).
$$

Assume that the polyhedron $X$ is represented as a union $X = \cup_\alpha X^\alpha$ of polyhedra. Assume further that $f^\alpha \in C^1(X^\alpha)$ for all $\alpha$, where $f^\alpha = f|_{X^\alpha}$, and also $f^\alpha_i = f^\beta_i$ on $X^\alpha \cap X^\beta$ for all $i, \alpha, \beta$. Then, clearly, $f \in C^1(X)$ with partial derivatives $f_i|_{X^\alpha} = f^\alpha_i$.

In Subsection \ref{subsec:SpaceOfCaps}, we identified the space $\C(D)$ with a convex polyhedron in $\R^\Sigma$, where the coordinates are the heights $(h_i)_{i \in \Sigma}$ of the cap. Since $h_i = 0$ for all $i \in \partial D$, we can consider $\C(D)$ as a polyhedron in $\R^{\Sigma \setminus \partial D}$.

\begin{prp} \label{prp:DerS}
The function $S$ is of class $C^2$ on $\C(D)$. Its partial derivatives are:
\begin{eqnarray}
\frac{\partial S}{\partial h_i} & = & \kappa_i, \label{eqn:DsDh}\\
\frac{\partial^2 S}{\partial h_i \partial h_j} & = & a_{ij}, \label{eqn:KiHj}\\
\frac{\partial^2 S}{\partial h_i^{\,2}} & = & - \sum_{j \in \Sigma} a_{ij}. \label{eqn:KiHi}
\end{eqnarray}
In (\ref{eqn:DsDh}) and (\ref{eqn:KiHj}) we assume $i, j \in \Sigma \setminus \partial D$ and $i \ne j$.
\end{prp}

\begin{proof}. First, let us show that $S$ is of class $C^2$ on every $\C^T(D)$ in the decomposition (\ref{eqn:Decomp}). Schl\"afli's formula \cite{Mil} for a Euclidean 3-polytope says that
$$
\sum_e \ell_e d\alpha_e = 0,
$$
where the sum extends over all edges of the polytope, $\ell_e$ is the length of the edge $e$, and $\alpha_e$ is the dihedral angle at $e$. By applying this to the prisms that constitute a generalized convex cap $C$, we obtain
$$
dS = \sum_{i \in \Sigma \setminus \partial D} \kappa_i dh_i
$$
which implies (\ref{eqn:DsDh}).

For simplicity, assume that the triangulation $T$ has no loops and no multiple edges. Note that the angle $\omega_i$ can be viewed as a function of the angles $\rho_{ij}, ij \in \E(T)$, as long as $h \in \C^T(D)$. Thus we have
\begin{eqnarray}
d \omega_i & = & \sum_{ij \in \E(T)} \frac{\partial \omega_i}{\partial \rho_{ij}} d \rho_{ij}, \label{eqn:dOmega}\\
d \rho_{ij} & = & \frac{\partial \rho_{ij}}{\partial h_i} dh_i + \frac{\partial \rho_{ij}}{\partial h_j} dh_j. \label{eqn:dRho}
\end{eqnarray}
From the trapezoid formed by the heights $h_i$ and $h_j$ it is easy to compute
\begin{equation} \label{eqn:DrhoDh}
\frac{\partial \rho_{ij}}{\partial h_i} = - \frac{1}{\ell_{ij} \sin \rho_{ij}} = -\frac{\partial \rho_{ij}}{\partial h_j}.
\end{equation}

Let $\omega_{ijk}$ be the dihedral angle at the $i$-th height in the prism $ijk$. Then we have
$$
\frac{\partial \omega_i}{\partial \rho_{ij}} = \frac{\partial \omega_{ijk}}{\partial \rho_{ij}} + \frac{\partial \omega_{ijl}}{\partial \rho_{ij}}.
$$
The spherical section of the prism $ijk$ at the vertex $i$ is a spherical triangle that has a side of length $\rho_{ij}$ with adjacent angles $\omega_{ijk}$ and $\alpha_{ij}$. With the help of the spherical Sine and Cosine theorems it is not hard to show
\begin{equation} \label{eqn:DomegaDrho}
\frac{\partial \omega_{ijk}}{\partial \rho_{ij}} = -\frac{\cot \alpha_{ij}}{\sin \rho_{ij}}.
\end{equation}
By substituting (\ref{eqn:DrhoDh}) and (\ref{eqn:DomegaDrho}) into (\ref{eqn:dRho}) and (\ref{eqn:dOmega}), we obtain
$$
d\omega_i = \sum_{ij \in \E(T)} \frac{\cot \alpha_{ij} + \cot \alpha_{ji}}{\ell_{ij} \sin^2 \rho_{ij}} (dh_i - dh_j) = \sum_{j \in \Sigma} a_{ij} (dh_i - dh_j).
$$
Formulas (\ref{eqn:KiHj}) and (\ref{eqn:KiHi}) follow from $\kappa_i = 2\pi - \omega_i$.

More generally, let $T$ contain loops and multiple edges. Let us denote by $e$ an oriented edge of $T$, by $a(e)$ and by $b(e)$ its initial and terminal vertex, respectively. Then formula (\ref{eqn:dOmega}) becomes
$$
d\omega_i = \sum_{a(e) = i} \frac{\partial \omega_i}{\partial \rho_e} d \rho_e.
$$
If $e$ is a loop, then $d \rho_e = 0$. Formulas (\ref{eqn:DrhoDh}) and (\ref{eqn:DomegaDrho}) with an obvious change of notation remain valid. Thus we again have $d \omega_i = \sum_{j \ne i} a_{ij}(dh_i - dh_j)$, where this time
$$
a_{ij} = \sum_e \frac{\cot \alpha_e + \cot \alpha_{-e}}{\ell_e \sin^2 \rho_e},
$$
whith the sum ranges over all edges joining $i$ and $j$.

Let a cap $C$ belong to several polyhedra $\C^T(D)$. Formulas (\ref{eqn:DsDh}), (\ref{eqn:KiHj}) and (\ref{eqn:KiHi}) yield the same result when computed for different geodesic triangulations associated with the cap $C$. Thus the partial derivatives up to order two of the function $S$ at $C$ are well-defined. This implies that the function $S$ is of class $C^2$ on the whole $\C(D)$.
\end{proof}

\subsection{Concavity of the total scalar curvature} \label{subsec:SConc}
\begin{thm} \label{thm:SConc}
The function $S$ is strictly concave.
\end{thm}
To prove this theorem, we use the information on the Hessian of $S$ obtained in the previous subsection.

\begin{dfn}
Let $C$ be a generalized convex cap with an associated triangulation $T$. The graph $\Gamma(C)$ is defined as follows. The vertex set of $\Gamma(C)$ is $\Sigma$; the edge set of $\Gamma(C)$ consists of the boundary edges of $T$ and of those interior edges $ij$ for which $\theta_{ij} < \pi$.
\end{dfn}
\begin{lem} \label{lem:Hess}
The Hessian $(\frac{\partial^2 S}{\partial h_i \partial h_j})$ is negatively semidefinite. The nullspace of $(\frac{\partial^2 S}{\partial h_i \partial h_j})$ is spanned by the vectors $v^K = (h_i^K)_{i \in \Sigma \setminus \partial D}$ defined as
$$
h_i^K = \left\{ \begin{array}{ll}
		1, & \mbox{for }i \in K;\\
		0, & \mbox{for }i \notin K,
                \end{array}
\right.
$$
where $K \subset \Sigma \setminus \partial D$ is a connected component of $\Gamma(C)$ disjoint with the boundary.
\end{lem}

\begin{proof}. Consider the quadratic form $\sum \frac{\partial^2 S}{\partial h_i \partial h_j} x_i x_j$. For simplicity, extend the summation over all $i, j \in \Sigma$ by putting $x_i = 0$ if $i \in \partial D$. By Proposition~\ref{prp:DerS} we have
\begin{eqnarray*}
\sum_{i,j \in \Sigma} \frac{\partial^2 S}{\partial h_i \partial h_j} x_i x_j & = &
\sum_{i \ne j} a_{ij} x_i x_j - \sum_{i \ne j} a_{ij} x_i^2 \\
& = & \sum_{i < j} (2a_{ij} x_i x_j - a_{ij} x_i^2 - a_{ij} x_j^2) \\
& = & - \sum_{i < j} a_{ij}(x_i - x_j)^2.
\end{eqnarray*}
Note that $a_{ij} \ge 0$, and $a_{ij} = 0$ if and only if in $\Gamma(C)$ there are no edges joining $i$ and $j$. Hence $\sum \frac{\partial^2 S}{\partial h_i \partial h_j} x_i x_j$ is always non-positive. It vanishes if and only if $x_i = x_j$ for $i$ and $j$ lying in one connected component of $\Gamma(C)$. Since $x_i = 0$ on the boundary component, the result follows. \end{proof}

\begin{proof} {\it of Theorem \ref{thm:SConc}.} By Lemma \ref{lem:Hess}, the function $S$ is concave. Assume that it is not strictly concave. Then there are two different generalized convex caps $C, C' \in \C(D)$ such that $S$ is linear on the segment joining $C$ and $C'$. Without loss of generality we can assume that $C, C' \in \C^T(D)$ for some geodesic triangulation $T$. This implies that the graph $\Gamma = \Gamma(C) \cup \Gamma(C')$ is embedded in $D$.

Let $h, h' \in \R^\Sigma$ be the height vectors of $C, C'$, respectively. Then the vector $(h'_i - h_i)_{i \in \Sigma \setminus \partial D}$ belongs to the nullspace of both $\frac{\partial^2 S}{\partial h_i \partial h_j} (C)$ and $\frac{\partial^2 S}{\partial h_i \partial h_j} (C')$. By Lemma \ref{lem:Hess}, this means that on every connected component of the graph $\Gamma$ the heights difference $h'_i - h_i$ is constant. If $\Gamma$ is connected, then we arrive to $C' = C$ since $h'_i = h_i = 0$ for $i \in \partial D$.

Let $\Gamma$ be disconnected. Graph $\Gamma$ cuts the disk $D$ in a number of pieces which we call the {\it faces} of $\Gamma$. A face may be adjacent to a single component of $\Gamma$ or to several ones. In the former case the face is necessarily simply connected, in the latter case it is multiply connected. Let us call two components of $\Gamma$ {\it neighboring} if they are adjacent to the same face. We show that the difference $h'_i - h_i$ is the same for any two neighboring components. This will imply $h' = h$ which contradicts to the assumption that $C'$ is different from $C$.

Let $K_1$ and $K_2$ be two neighboring components of $\Gamma$, and let $F$ be the face adjacent to both $K_1$ and $K_2$. Consider the PL extensions $\widetilde{h}$ and $\widetilde{h'}$ of $h$ and $h'$, respectively (see Lemma \ref{lem:CapsAreFunctions}). Function $\widetilde{h}$ is linear on the face $F$. Thus the gradient of $\widetilde{h}$ defines a parallel vector field on $F$. Since $F$ is multiply connected, it contains a simple loop that is non-contractible in $D \setminus \Sigma$. The parallel transport along this loop is the rotation by the angle equal to the sum of the curvatures of the singularities enclosed by the loop. This sum is less than $2\pi$, because of $\Sigma \cap \partial D \ne \emptyset$. It follows that $F$ cannot carry a nonzero parallel vector field. Thus $\widetilde{h}$ is constant on $F$. By the same reason, $\widetilde{h'}$ is constant on $F$. Therefore, $h' - h$ is constant on the vertices of $F$. Since $F$ has a vertex in $K_1$ and a vertex in $K_2$, we conclude that $h' - h$ is constant on the vertices of $K_1 \cup K_2$. \end{proof}

\noindent\textbf{Remark} Consider the generalized cap with all heights zero. Then $a_{ij} = 0$ for all $i,j$ and the Hessian vanishes identically. If we put $h_i = \epsilon$ for all $i \notin \partial D$, where $\epsilon \le d_i$ for all $i$, then $h \in \C(D)$ and thus defines a cap $C_\epsilon$. It is easy to find a polyhedral disk $D$ such that the graph $\Gamma(C_\epsilon)$ is disconnected. This example shows that the points of degeneration of the Hessian might be non-isolated.

\subsection{Applications to infinitesimal and global rigidity}
Let $n = |\Sigma \setminus \partial D|$ be the number of interior singularities of disk $D$. There is a map
$$
\begin{array}{rcl}
K_D: \C(D) & \to & \R^n,\\
r & \mapsto & \kappa,
\end{array}
$$
where $r = (r_1,\ldots, r_n)$ and $\kappa = (\kappa_1, \ldots, \kappa_n)$. Different kinds of rigidity for caps with the given boundary are defined by the extent to which a cap is determined by its curvatures, that is whether $K_D$ is locally or globally injective or not.

\begin{dfn}
One says that generalized convex caps with the upper boundary $D$ are {\rm globally rigid} iff the map $K_D$ is injective.
\end{dfn}

\begin{thm} \label{thm:GenRig}
For any convex polyhedral disk $D$, generalized convex caps with the upper boundary $D$ are globally rigid.
\end{thm}
\begin{proof}. By \eqref{eqn:DsDh},
$$
K_D = \grad S,
$$
and by Theorem \ref{thm:SConc} the function $S$ is strictly concave. Its domain $\C(D)$ is compact and convex by Theorem \ref{thm:C(D)}. Thus Proposition \ref{prp:ConvHomeo} below implies that the map $K_D$ is a homeomorphism onto the image. \end{proof}

\begin{prp} \label{prp:ConvHomeo}
Let $f \in C^1(X)$ be a strictly convex or strictly concave function on a compact convex set $X \subset \R^n$. Then the map $\grad f: X \to \R^n$ is a homeomorphism onto the image.
\end{prp}

The same technique was recently used by Feng Luo in \cite{Luo06} to prove a collection of rigidity theorems in dimension 2.

\begin{proof}. Since $X$ is compact and $\grad f$ is continuous, it suffices to show that $\grad f$ is injective. Let $x$ and $y$ be two different points in $X$. The restriction of $f$ to the segment $xy$ is a convex $C^1$-function. Thus the directional derivative
$$
\frac{\partial f}{\partial \xi} (\lambda x + (1-\lambda) y),
$$
where $\xi = (y-x)/\|y-x\|$, is a monotone function of $\lambda \in [0,1]$. In particular,
$$
\frac{\partial f}{\partial \xi}(x) \ne \frac{\partial f}{\partial \xi}(y).
$$
Since $\frac{\partial f}{\partial \xi} = \langle \grad f, \xi \rangle$, it follows that $\grad f(x) \ne \grad f(y)$. \end{proof}

\begin{dfn} \label{dfn:InfRigGen}
Let $C \in \C(D)$ be a generalized convex cap. One says that $C$ is {\rm infinitesimally rigid} iff the Jacobian of the map $K_D$ at $C$ has full rank.
\end{dfn}

Since the Jacobian of $K_D$ equals the Hessian of $S$, Lemma \ref{lem:Hess} provides a necessary and sufficient condition for $C$ to be infinitesimally rigid. Namely, it is the case iff the graph $\Gamma(C)$ is connected. The remark at the end of Subsection \ref{subsec:SConc} provides examples of generalized convex caps which are not infinitesimally rigid. On the other hand, the argument in the proof of Theorem \ref{thm:SConc} shows that every inifinitesimally flexible cap has a multiply connected face, and this face is parallel to the base of the cap. It is not hard to show that infinitesimally flexible caps lie on the boundary of the space $\C(D)$.

Now let us turn to classical convex caps and convex polytopes. If $C$ is a classical convex cap with the upper boundary $D$, then the graph $\Gamma(C)$ is clearly connected. Thus $C$ is infinitesimally rigid in the sense of Definition \ref{dfn:InfRigGen}. However, usually one defines the infinitesimal rigidity for polyhedral objects in a different way.

\begin{dfn} \label{dfn:InfRigCla}
Let $C$ be a classical convex cap without vertical faces. Triangulate its faces by diagonals. An {\rm infinitesimal flex} of $C$ is a collection of vectors $p'_1, \ldots, p'_m$ such that
\begin{enumerate}
\item $(p'_i - p'_j) \cdot (p_i - p_j) = 0$ if $ij$ is an edge of $C$ or one of the added diagonals, where $\cdot$ denotes the scalar multiplication;
\item $p_i$ belongs to the horizontal plane if $i$ is a boundary singularity.
\end{enumerate}
The cap $C$ is called {\rm infinitesimally rigid} if for any triangulation of its faces every infinitesimal flex is a restriction of an infinitesimal rigid motion in~$\R^3$.
\end{dfn}

\begin{lem} \label{lem:InfRigEq}
Let $C$ be a classical convex cap without vertical faces and with the upper boundary $D$. Then $C$ is infinitesimally rigid as a generalized convex cap, see Definition \ref{dfn:InfRigGen}, iff $C$ is infinitesimally rigid in the sense of Definition \ref{dfn:InfRigCla}.
\end{lem}
\begin{proof}. Let $(p'_i)$ be an infinitesimal flex  of $C$ for some triangulation of the faces. Let $h'_i$ be the vertical component of the vector $p'_i$. If $(p'_i)$ is not a restriction of a rigid motion, then it can be shown that $h' = (h'_1, \ldots, h'_n)$ is a non-zero vector in $\R^{\Sigma \setminus \partial D}$. But clearly, $h' \in \ker \left(\frac{\partial \kappa_i}{\partial h_j}\right)$. Thus the classical flexibility implies the generalized flexibility.

Assume $h' \in \ker \left(\frac{\partial \kappa_i}{\partial h_j}\right)$ is a non-zero vector. Choose a triangulation of the faces of $C$. We claim that there is an infinitesimal flex $(p'_i)$ that has $(h'_i)$ as its vertical components. Note that the horizontal components for the endpoints of an edge determine the horizontal component for the vertex opposite to this edge. This allows to determine the horizontal components for all vertices starting with an arbitrary edge. Despite many possible ways to arrive at the vertex $i$ starting from the given edge, the vector $p'_i$ is well-defined due to $h' \in \ker \left(\frac{\partial \kappa_i}{\partial h_j}\right)$. Thus the generalized flexibility implies the classical flexibility. \end{proof}

\begin{thm} \label{thm:InfRigCap}
Classical convex caps that have dimension 3 are infinitesimally rigid in the sense of Definition \ref{dfn:InfRigCla}.
\end{thm}
If the cap is a polygon orthogonal to the base plane, then this polygon can be bent along a diagonal. If we want the infinitesimal rigidity to hold in this case, Definition \ref{dfn:InfRigCla} must be appropriately modified.
\begin{proof}. Let $C$ be a classical convex cap. The graph $\Gamma(C)$ is connected, therefore $C$ is infinitesimally rigid as a generalized cap. By Lemma \ref{lem:InfRigEq}, if $C$ has no vertical faces, then it is infinitesimally rigid in the sense of Definition~\ref{dfn:InfRigCla}.

Let $C$ be a cap that has vertical faces and is 3-dimensional. Denote by $D^\wedge$ the upper boundary of $C$ minus the vertical faces. For any infinitesimal flex of $C$ the vertical components of the vectors $p'_i$ for $i \in \partial D^\wedge$ vanish. Let $a_i$ be the height of $i$ in $C$ for $i \in \partial D^\wedge$. Consider the intersection of the space $\C(D)$ with the affine plane $\{h_i = a_i \mbox{ for all } i \in \partial D^\wedge\}$, and consider only the coordinates $\{\kappa_i \,|\, i \notin \partial D^\wedge\}$ of the map $K_D$. Since the Jacobian of $K_D$ at $C$ is a negative definite symmetric matrix, the Jacobian of the restricted $K_D$ is also negative definite, in particular non-degenerate. Hence by the first part of the proof of Lemma \ref{lem:InfRigEq} the cap $C$ is infinitesimally rigid in the sense of Definition \ref{dfn:InfRigCla}. \end{proof}

For polytopes in $\R^3$ with triangular or triangulated faces, infinitesimal flexes are defined exactly as in Definition \ref{dfn:InfRigCla}, except that now we have no condition~2. Infinitesimal rigidity of convex polytopes was proved by Max Dehn in \cite{Dehn}. We present here a new proof based on Theorem \ref{thm:InfRigCap} and on the projective invariance of the infinitesimal rigidity.

\begin{thm}[Dehn] \label{thm:Dehn}
Convex 3-polytopes are infinitesimally rigid.
\end{thm}
\begin{proof}. Let $P \subset \R^3$ be a convex polytope. It is an old result \cite{Dar93}, \cite{Sau35} that for smooth surfaces the infinitesimal rigidity is a projectively invariant property. The arguments of Darboux and Sauer extend to the case of polyhedral surfaces.

So let us consider a projective transformation $\phi$ of $\R^3$ that maps a vertex $p_0$ of $P$ to the point with the homogeneous coordinates $[0,0,1,0]$, where $\R^3$ is identified with the affine part of $\R P^3$ via $(x,y,z) = [x,y,z,1]$. Then the faces of $P$ that contain $p_0$ become vertical. The map $\phi$ can be chosen so that $\phi(p_i) \in \R_+^3$ for all vertices of $P$ except $p_0$ and $C = \phi(P) \cap \R_+^3$ is a convex cap. By Lemma 5.1 in \cite{Scl07}, with every infinitesimal flex of $P$ that fixes $p_0$ one can associate an infinitesimal flex of $C$. Besides, the trivial flexes of $P$ go exactly to the trivial flexes of $C$. Thus Theorem \ref{thm:InfRigCap} implies that $P$ is infinitesimally rigid. \end{proof}

\noindent\textbf{Remarks} This proof of Theorem \ref{thm:Dehn} is due to Jean-Marc Schlenker. In \cite{Scl07}, he goes further by proving that the Hessian of $S$ is negative definite also for weakly convex caps, that is caps that are not necessarily convex, but whose vertices are extremal points of a convex polytope. As in this Subsection, this implies that weakly convex caps are infinitesimally rigid. From the Darboux-Sauer theorem Schlenker deduces that weakly convex polytopes that are star-shaped with respect to a vertex are infinitesimally rigid.

In \cite{BI}, a proof of Dehn's theorem was given that uses the total scalar curvature $S$ for generalized convex polytopes. There the function $S$ is not concave, but has one positive direction. The rigidity follows from the fact that the corank of the Hessian of $S$ equals 3, that is the dimension of the space of infinitesimal rigid motions. The signature of $S$ is computed using the connection of the total scalar curvature to the volume of the dual as established in \cite{BI} and from the Alexandrov-Fenchel inequalities.

\section{Proof of Theorem \ref{thm:Alex}} \label{sec:Proof}
Here is the idea of the proof. We have a strictly concave function $S$ on the space $\C(D)$ of convex generalized caps. A cap $C$ that lies in the interior of $\C(D)$ is a classical convex cap if and only if it has curvatures $\kappa_i = 0$ for all $i$. On the other hand, we have $\grad S = \kappa$. Thus, in the interior of $\C(D)$ the critical points of $S$ correspond to the classical convex caps. Due to the concavity of the function $S$, any interior critical point is a point of local maximum of $S$. And due to the convexity of $\C(D)$, a local maximum is the global maximum and is therefore unique.

The problem is that $S$ may attain its maximum at a boundary point of $\C(D)$. Thus we need to characterize those boundary points that may maximize $S$ as well as those that may produce a classical cap. Subsection \ref{subsec:Proof} deals mostly with this problem. In Subsection \ref{subsec:Volkov} we give an account of Volkov's proof.

\subsection{Classical caps maximize the total scalar curvature} \label{subsec:Proof}
\begin{dfn} \label{dfn:DegPart}
Let $C$ be a generalized convex cap. The {\rm degenerate part} $C^|$ of $C$ is defined as the union of the degenerate prisms of $C$. The {\rm non-degenerate part} $C^\wedge$ of $C$ is the union of the non-degenerate prisms.
\end{dfn}

\begin{lem} \label{lem:ClassCurv}
A generalized convex cap $C$ is a classical convex cap if and only if its curvatures satisfy the following conditions:
\begin{eqnarray}
\kappa_i = 0 & \mbox{if} & h_i < d_i; \label{eqn:h<d}\\
\kappa_i \ge 0 & \mbox{if} & h_i = d_i. \label{eqn:h=d}
\end{eqnarray}
Here $d_i$ is the distance from the singularity $i$ to $\partial D$ in the metric of $D$.
\end{lem}
\begin{proof}. For a generalized convex cap $C$, consider the intersection $D^| = C^| \cap D$. It consists of the upper bases of the degenerate prisms of $C$. Equivalently, $D^|$ is the closure of the set of $x \in D$ such that $|\grad \widetilde{h}(x)| = 1\}.$ Therefore if $x \in D^|$, then also $y \in D^|$ for all $y$ on the shortest geodesic joining $x$ with $\partial D$. It follows that $D^|$ can be represented as a union of polygons
\begin{equation} \label{eqn:Walls}
D^| = \bigcup_s P_s,
\end{equation}
where every $P_s$ has exactly one side $L_s$ in $\partial D$, and the orthogonal projection to $L_s$ maps the upper ridge $\partial P_s \setminus L_s$ injectively into $L_s$. See Figure \ref{fig:Wall}.

\begin{figure}[ht]
\begin{center}
\begin{picture}(0,0)%
\includegraphics{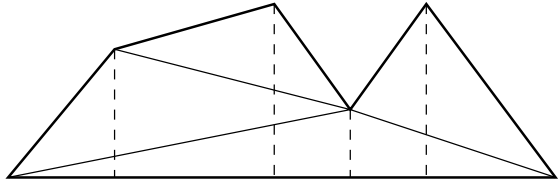}%
\end{picture}%
\setlength{\unitlength}{4144sp}%
\begingroup\makeatletter\ifx\SetFigFont\undefined%
\gdef\SetFigFont#1#2#3#4#5{%
  \reset@font\fontsize{#1}{#2pt}%
  \fontfamily{#3}\fontseries{#4}\fontshape{#5}%
  \selectfont}%
\fi\endgroup%
\begin{picture}(3284,1079)(429,-1208)
\end{picture}%
\end{center}
\caption{A polygon $P_s$ from the decomposition (\ref{eqn:Walls}). Corresponding degenerate prisms are also shown.}
\label{fig:Wall}
\end{figure}

Clearly, every $P_s$ is isometric to the union of the degenerate prisms whose upper bases lie in $P_s$. These prisms contribute $\pi$ to the angles $\omega_i$ around the vertical edges from the singularities $i$ in the upper ridge $\partial P_s \setminus L_s$.

Note that the polygons in (\ref{eqn:Walls}) can have points in common, moreover a polygon can have identifications on the boundary, that is it may be immersed and not embedded in $D$. However, if the curvatures of $C$ satisfy the condition (\ref{eqn:h=d}), then there are only the following three possibilities:

1) Every polygon is embedded, and no two polygons have a point in common, except maybe points in $\partial D$.

2) In (\ref{eqn:Walls}) there are only two polygons, they are equal, and $D$ is obtained by identifying their upper ridges. In particular, all prisms are degenerate.

3) There is only one polygon $P$, it has a vertical axis of symmetry, and $D$ is obtained from $P$ by identifying the corresponding points of the upper ridge. Again, all prisms are degenerate.

Now, let the curvatures of $C$ satisfy both conditions (\ref{eqn:h<d}) and (\ref{eqn:h=d}). If $C^\wedge \ne \emptyset$, then we have the situation 1). Hence $C^\wedge$ is isometric to a  polytope in $\R^3_+$ that projects to its lower base and is convex except maybe at the vertical edges. The polygons $P_s$ are glued to the vertical faces of the polytope. Since they contribute $\pi$ to the angles at the vertical edges, the polytope is a convex cap.

If we have the situation 2) or 3), then $C$ is isometric to a 2-dimensional classical cap, where in 3) the cap has a vertical edge.

Thus conditions (\ref{eqn:h<d}) and (\ref{eqn:h=d}) are sufficient for a cap $C$ to be classical. Proof of the necessity is not hard, and we omit it. \end{proof}


\begin{proof} {\it of Theorem \ref{thm:Alex}.}

{\it Existence.} Let $C \in \C(D)$ be a maximum point of the function $S$. We claim that $C$ is a classical convex cap.

If $C$ lies in the interior of $\C(D)$, then we have $\kappa(C) = \grad S(C) = 0$. Therefore $C$ is a classical convex cap with no vertical faces. If $C$ lies on the boundary of $\C(D)$, then we need to show that $C$ satisfies the properties (\ref{eqn:h<d}) and (\ref{eqn:h=d}).

Assume that (\ref{eqn:h<d}) or (\ref{eqn:h=d}) is violated. First, consider the case when there is an $i$ such that $\kappa_i < 0$. Let us show that the height $h_i$ can be decreased by a small amount so that $h$ does not leave the set $\C(D)$. Since $\frac{\partial S}{\partial h_i} = \kappa_i < 0$, the function $S$ grows during this deformation, and this contradicts the choice of the cap $C$. We have to show that the inequalities (\ref{eqn:ConcQuadr}), (\ref{eqn:CloGeod}), and (\ref{eqn:hInt}) remain valid when $h_i$ decreases. For (\ref{eqn:hInt}) it is obvious. A small decrease of $i$ violates the inequality (\ref{eqn:ConcQuadr}) only if the equality $h_i = \ext_{jkl}(i)$ holds at the beginning. This means that $\theta_{ij} = \pi$, and thus there is a face with the angle at $i$ greater or equal $\pi$. But this contradicts Lemma \ref{lem:ConvFaces}. Similarly, if we violate (\ref{eqn:CloGeod}), then at the beginning we have an equality there. But clearly, $h_i = h_j$ implies that the curvature $\kappa_i$ is equal to the angle defect of the singularity $i$ in the metric of $D$. This contradicts the assumption $\kappa_i < 0$. Finally, a decrease of $h_i$ does not violate the inequalities of type (\ref{eqn:ConcQuadr}) and (\ref{eqn:CloGeod}), if the singularity $i$ plays the role of $j$, $k$ or $l$ there.

Thus we may assume that $\kappa_i \ge 0$ for all $i$, and that there is an $i$ such that $\kappa_i > 0$ and $h_i < d_i$. Increase by the same small amount all of those $h_i$ that are smaller than $d_i$. If $h$ does not leave $\C(D)$, then this deformation increases $S$ that contradicts the choice of the cap $C$. So assume that one of the inequalities (\ref{eqn:ConcQuadr}), (\ref{eqn:CloGeod}), and (\ref{eqn:hInt}) gets violated. The inequality (\ref{eqn:hInt}) remains valid, if the deformation is sufficiently small. The inequality (\ref{eqn:CloGeod}) gets violated only if it is an equality at the beginning, and $h_j$ gets increased but $h_i$ remains constant. The equality $h_i = h_j$ implies that $\widetilde{h}$ is constant on the triangle $jij$. Thus $i$ lies in the non-degenerate part of $D$, and we have $h_i < d_i$. This contradicts to the assumption that $h_i$ remains constant during the deformation. Now let (\ref{eqn:ConcQuadr}) get violated. Then at the initial moment we have $h_i = \ext_{jkl}(i)$, and $h_i$ remains constant during the deformation, but some of the other heights grow. Since $h_i$ remains constant, we have $h_i = d_i$. Thus $i$ lies in the degenerate part $D^|$ of $D$. On the other hand, the quadrilateral $ikjl$ is not contained in $D^|$, otherwise all of the vertices $j$, $k$, $l$ lie in the degenerate part, and the corresponding heights remain constant during the deformation. Thus both the degenerate and non-degenerate prisms contribute $\pi$ to the angle $\omega_i$. So we have $\kappa_i \le 0$. Since the face containing the quadrilateral $ikjl$ has the angle at least $\pi$ at $i$, this contradicts Lemma \ref{lem:ConvFaces}.

{\it Uniqueness.} Let $C \in \C(D)$ be isometric to a classical convex cap. Let us show that $C$ is a point of local maximum for the function $S$. The curvatures of $C$ satisfy conditions (\ref{eqn:h<d}) and (\ref{eqn:h=d}). If $h_i < d_i$ for all interior singularities $i$, then $\kappa_i = 0$ implies that $C$ is a critical point of $S$, and thus a point of local maximum due to the concavity of $S$. If $h_i = d_i$ for some $i$, then we have $\frac{\partial S}{\partial h_i} = \kappa_i \ge 0$. On the other hand, $h_i \le d_i$ on $\C(D)$. Therefore directional derivatives of $S$ at $C$ in directions pointing inside $\C(D)$ are non-positive. Again, the concavity of $S$ implies that $C$ is a point of local maximum.

By Theorem \ref{thm:SConc}, the function $S$ is strictly concave. Therefore it has only one local maximum on the convex space $\C(D)$, and the uniqueness follows.

Actually, instead of Theorem \ref{thm:SConc} a more simple Lemma~\ref{lem:Hess} can be applied. Assume that $C_1, C_2 \in \C(D)$ are two different caps isometric to classical ones. In a classical convex cap the graph $\Gamma(C_i)$ is connected. Thus the Hessian of $S$ is positively defined at both $C_1$ and at $C_2$. Join the points $C_1$ and $C_2$ by a straight segment in $\C(D) \subset \R^{\Sigma \setminus \partial D}$. Since $C_1$ and $C_2$ are points of local maximum for $S$, the function $S$ must be constant on the segment $C_1C_2$. But this contradicts to the positivity of the second derivative at the segment endpoints. \end{proof}


Now we prove the lemma used in the proof of the existence of a convex cap. Recall that we defined faces of $C$ as the connected components of $D \setminus \Gamma(C)$, where $\Gamma(C)$ is the graph obtained from the skeleton of an associated triangulation by removing all of the interior edges $ij$ with $\theta_{ij} = \pi$. In general, a face may be non-convex.
\begin{lem} \label{lem:ConvFaces}
Let $C$ be a generalized convex cap. Assume that $\kappa_i \le 0$ for some interior singularity $i$. Then all of the faces of $C$ incident to the vertex $i$ have angles at $i$ less than $\pi$.
\end{lem}
\begin{proof}. Consider the spherical section $C_i$ of the cap $C$ at the vertex $i$. It is glued from spherical triangles with a common vertex $O$. The total angle at $O$ equals $2\pi - \kappa_i \ge 2\pi$. Maximal geodesic arcs on the boundary of $C_i$ correspond to the face angles of $C$ at $i$. Assume that there is a geodesic arc $AB \subset \partial C_i$ of length $\pi$. Then the arcs $OA$ and $OB$ bound a spherical lune in $C_i$, see the left part of the Figure \ref{fig:SpherSec}.

\begin{figure}[ht]
\begin{center}
\begin{picture}(0,0)%
\includegraphics{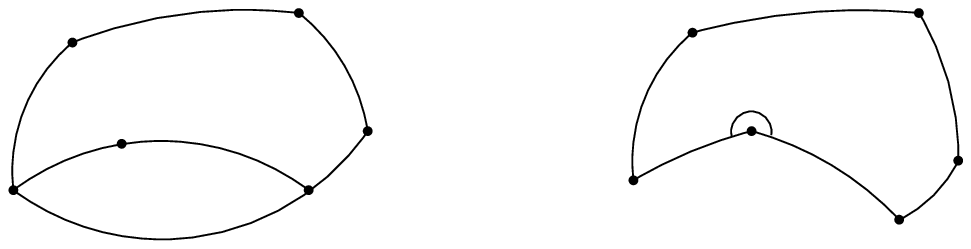}%
\end{picture}%
\setlength{\unitlength}{4144sp}%
\begingroup\makeatletter\ifx\SetFigFont\undefined%
\gdef\SetFigFont#1#2#3#4#5{%
  \reset@font\fontsize{#1}{#2pt}%
  \fontfamily{#3}\fontseries{#4}\fontshape{#5}%
  \selectfont}%
\fi\endgroup%
\begin{picture}(4582,1145)(271,-596)
\put(271,-466){\makebox(0,0)[lb]{\smash{{\SetFigFont{12}{14.4}{\rmdefault}{\mddefault}{\updefault}{\color[rgb]{0,0,0}$A$}%
}}}}
\put(1801,-466){\makebox(0,0)[lb]{\smash{{\SetFigFont{12}{14.4}{\rmdefault}{\mddefault}{\updefault}{\color[rgb]{0,0,0}$B$}%
}}}}
\put(3106,-421){\makebox(0,0)[lb]{\smash{{\SetFigFont{12}{14.4}{\rmdefault}{\mddefault}{\updefault}{\color[rgb]{0,0,0}$A$}%
}}}}
\put(4546,-556){\makebox(0,0)[lb]{\smash{{\SetFigFont{12}{14.4}{\rmdefault}{\mddefault}{\updefault}{\color[rgb]{0,0,0}$B$}%
}}}}
\put(901,-286){\makebox(0,0)[lb]{\smash{{\SetFigFont{12}{14.4}{\rmdefault}{\mddefault}{\updefault}{\color[rgb]{0,0,0}$O$}%
}}}}
\put(3556,119){\makebox(0,0)[lb]{\smash{{\SetFigFont{12}{14.4}{\rmdefault}{\mddefault}{\updefault}{\color[rgb]{0,0,0}$\pi-\kappa_i$}%
}}}}
\put(3736,-241){\makebox(0,0)[lb]{\smash{{\SetFigFont{12}{14.4}{\rmdefault}{\mddefault}{\updefault}{\color[rgb]{0,0,0}$O$}%
}}}}
\end{picture}%
\end{center}
\caption{The spherical section of a cap at a vertex, and its development after cutting out the lune $OAB$.}
\label{fig:SpherSec}
\end{figure}

Denote by $BA$ the complement to the arc $AB$ in the boundary of $C_i$. Remove from $C_i$ the spherical lune $OAB$ and develop the rest onto the sphere. $BA$ becomes a piecewise geodesic of positive turn. Due to $\kappa_i \le 0$, the angle under which $BA$ is seen from $O$ is greater or equal $\pi$ (note that it may be even greater than $2\pi$). This implies $|BA| \ge |OA| + |OB| = \pi$. Since we have $|AB| = \pi$, the boundary of $C_i$ has the total length at least $2\pi$. But this length equals the total angle around the singularity $i$ in the metric of $D$. This cotradicts the assumption that $D$ is a convex polyhedral disk. \end{proof}

\subsection{Volkov's proof} \label{subsec:Volkov}
Here we give a short account of Volkov's proof \cite{Vol60} of the existence part of Theorem \ref{thm:Alex}. Its main ingredients are incorporated in our proof.

Volkov considers the set of all generalized convex caps with non-negative curvatures $\kappa_i$, see Definitions \ref{dfn:GCC} and \ref{dfn:kappa}. He takes the cap that maximizes the sum $\sum_i h_i$ of the heights and proves that it has zero curvatures in its degenerate part, see Definition \ref{dfn:DegPart}. Hence it is a classical convex cap. The claim that for the cap with the maximum sum of heights all curvatures vanish is proved by contradiction: let $C$ be a generalized convex cap that maximizes $\sum_i h_i$ and assume that $\kappa_i > 0$ for some $i$ in the non-degenerate part. Then one increases by a small amount all of the heights $h_i$ for which $\kappa_i > 0$ and $i$ is not in the non-degenerate part, and obtains a generalized convex cup of non-positive curvature with a greater sum of heights. This is exactly what we do in the second part of our existence proof. To show that the heights with positive curvatures may be increased, Volkov proves Lemma \ref{lem:ConvFaces} and gives a subtle geometric argument that shows how the upper boundary of the cap must be retriangulated.

\end{document}